\documentclass[journal]{IEEEtran}

\usepackage{cite}

\ifCLASSINFOpdf
\usepackage[pdftex]{graphicx}
\graphicspath{{../pdf/}{../jpeg/}}
\DeclareGraphicsExtensions{.pdf,.jpeg,.png}
\else
\fi

\usepackage[cmex10]{amsmath}

\usepackage{amssymb}

\usepackage{comment}
\usepackage{tikz}

\usepackage{accents}

\usepackage{algorithm}
\usepackage{algorithmic}

\usepackage{mathtools}

\begin{document}
\title{Dynamic Polytopic Template Approach to Robust Transient Stability Assessment}
\author{
Dongchan~Lee,~\IEEEmembership{Student Member,~IEEE,}	Konstantin~Turitsyn,~\IEEEmembership{Member,~IEEE}
\thanks{D. Lee and K. Turitsyn are with the Department of Mechanical Engineering, Massachusetts Institute of Technology, Cambridge, MA 02139, USA (email: dclee@mit.edu; turitsyn@mit.edu).}
\thanks{This work was supported by the NSF  awards  1554171  and  1550015 and Advanced Grid Modeling Program of the Office of Electricity within the U.S Department of Energy.}
}

\markboth{Submitted for publication. This version: May 2, 2017}{}

%



\maketitle

\begin{abstract}
Transient stability assessment of power systems needs to account for increased risk from uncertainties due to the integration of renewables and distributed generators. The uncertain operating condition of the power grid hinders reliable assessment of transient stability. Conventional approaches such as time-domain simulations and direct energy methods are computationally expensive to take account of uncertainties. This paper proposes a reachability analysis approach that computes bounds of the possible trajectories from uncertain initial conditions. The eigenvalue decomposition is used to construct a polytopic template with a scalable number of hyperplanes that is guaranteed to converge near the equilibrium. The proposed algorithm bounds the possible states at a given time with a polytopic template and solves the evolution of the polytope over time. The problem is solved with linear programming relaxation based on outer-approximations of nonlinear functions, which is scalable for large scale systems. We demonstrate our method on IEEE test cases to certify the stability and bound the state trajectories.
\end{abstract}


\IEEEpeerreviewmaketitle

\section{Introduction}

Transient stability assessment is one of the most critical tools for ensuring the reliability of power systems against any potential disturbances and contingencies \cite{kundur94,machowski11}. Due to the nonlinearity of system dynamics, the interconnected generators may not be able to resynchronize their frequencies. The loss of synchronization may lead to outages and blackouts, incurring high economic and societal cost. The system operators derive the operational limits to ensure the security of the system under N-1 contingency scenarios. These limits are derived in an off-line setting prior to the operation due to the size and complexity of the system \cite{kundur00, lee15}.

To minimize the generation cost, the operation takes place near the system limit where uncertainties can play a critical role. On the other hand, the rapid integration of renewables and power electronic devices has increased uncertainties in the system, which in turn increases the possibility of violating security constraints \cite{milano13}. However, taking account of uncertainties significantly increases the computational cost for large-scale power systems using conventional methods. This forces the operational limit to rely on the safety margin, which forces the system to operate in a conservative manner. Therefore, it is important for system operators to be equipped with an accurate and efficient assessment tool that can take account of the growing concerns with uncertainties. Our paper aims to provide a novel technique that is tractable and certifiable for robust transient stability assessment under the uncertain operating condition.

The most wide-spread approach for transient stability assessment has been time-domain simulations that test possible contingencies in the N-1 security set at various operating conditions \cite{pavella12}. The time domain simulation gives the solution to the trajectory, given that we know accurate initial condition and system parameters. However, the incorporation of uncertainties requires sampling-based approaches, such as the Monte-Carlo technique \cite{hockenberry04}. An alternative approach is based on the direct energy method, which uses the energy function to certify stability \cite{chiang15}. This approach was generalized to the Lyapunov Functions Family with systematic derivation of the stable region \cite{vu16_lyap,vu16_robust}. However, this approach often produces a conservative boundary, which may not cover the region of interest. In addition, the uncertainties in power injection for the transient stability assessment have received much attention. \cite{dong12,dhople13}.

Recently, a reachability assessment approach has been proposed in \cite{sassi12, dang12} for polynomial systems. We share the underlying idea on using a template polytope to bound the discretized time steps. However, the algorithm in the paper is scalable only to relatively low-order polynomial systems. We propose a more tractable and scalable algorithm for transient stability assessment. A similar idea, known as interval analysis, was used to bound the measurement and numerical error \cite{moore09, jaulin02}. This approach was used for power systems to bound the trajectories under disturbances \cite{pico13}. However, the convergence of the bounded states in interval analysis has not been possible because it is usually conservative. The intervals form a box that typically does not converge in under-damped systems. To alleviate this major limitation, we introduce a novel approach that constructs the bounding template using a polytope, which incorporates the geometrical characterization of system dynamics. Our construction is also strongly related to the contraction theory \cite{lohmiller98, singh17} but cannot be naturally recovered from the traditional metrics and norms discussed in the literature.

The reachability analysis approach introduced in this paper is based on solving linear programming and can be applied to high-order generator models in power systems represented as differential-algebraic equations. Our approach is apply discretization to the original differential equations and use a polytopic template to iteratively approximate the reachability region at every time step. The resulting sequence of polytopes is guaranteed to contain all possible trajectories, while the contraction of the polytope is used to certify stability. An explicit exit condition of our algorithm is derived by computing the largest invariant template polytope that is uniformly bounded. The resulting reachability set can be naturally used to certify that the system never violates any operational constraints, like frequency or angle limits, during transient events. Since our algorithm is based on linear programming, we can solve very large problems efficiently \cite{bertsimas97}. We conduct a study on the third-order synchronous generator model and demonstrate the effectiveness of our solution.

Our paper is organized as follows. In Section II, we formulate reachability analysis and give instructions on constructing the polytopic template. Section III presents the outer-approximation to solve the optimization problem formulated in Section II. Section IV shows an illustrative example on a 2 bus system as well as case studies on IEEE 14 and 39 bus systems with third-order generator model. We continue our discussions in Section V  and conclude in Section VI with future directions.

\section{Problem Formulation and Preliminaries}

\begin{figure}[!htbp]
	\centering
	\includegraphics[width=3.2in]{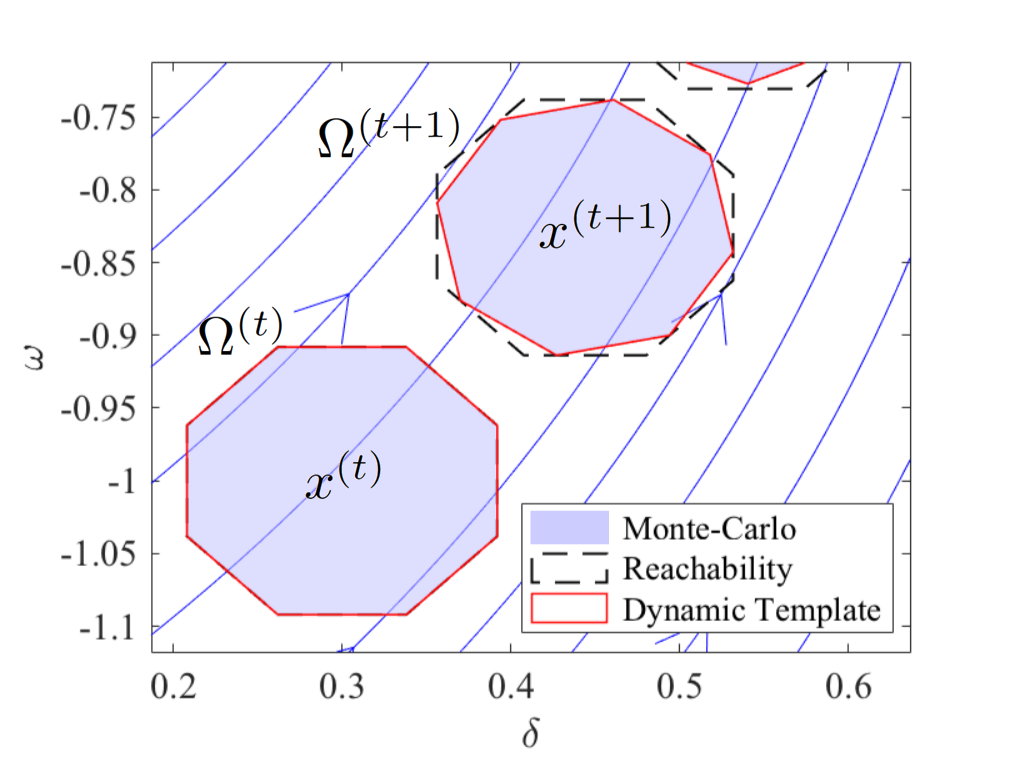}
	\caption{The reachability analysis is illustrated and compared with the Monte-Carlo method. The region for the Monte-Carlo simulation was computed by taking the convex hull of simulated points at $t+1$.}
	\label{fig_prob_form}
\end{figure}

Power system dynamics can be generally written with differential algebraic equations (DAE) in the following form:
\begin{equation}
	\begin{aligned}
		\dot{x}&=f(x,y) \\
		0&=g(x,y)
	\end{aligned}
\end{equation}
where $x$ and $y$ are the variables for the differential equation and algebraic equations, respectively. The generator-related variables such as rotor angle and rotor frequency are often the differential variables while the current flows over transmission lines and voltage at buses are the algebraic variables. Our formulation uses the explicit Euler's method with time step $\Delta t$ to discretize the system:
\begin{equation}
	\begin{aligned}
		x^{(t+1)}&=x^{(t)} +\Delta t f(x^{(t)},y^{(t)}) \\
		0&=g(x^{(t)},y^{(t)}).
	\end{aligned}
	\label{eqn_forward_Euler}
\end{equation}
This discretization turns the continuous ordinary differential equation into a set of algebraic equations, and we will exploit this method to develop our algorithm. The polytope will be used to bound the possible operating points over the differential variables, and we denote this set by
\begin{equation}
	\Omega^{(t)} = \{x^{(t)} \ | \ A^{(t)}x^{(t)} \leq b^{(t)}\}
\end{equation}
where $A^{(t)} \in R^{n\times m}$ and $b^{(t)} \in R^{n}$ define the polytope at time $t$. This form is an extension of simple interval bounds over each variable. This polytope generalizes intervals and capture the relationship between generator states as linear constraints. The reachability analysis aims to find this bound of all possible states at a given time, and numerically computes that the set, $\Omega^{(t)}$, converges towards the equilibrium. Since our set contains all possible states, the convergence of this set certifies the stability of every initial condition that was in the initial polytope. This idea is illustrated in Figure \ref{fig_prob_form}, where the octagon at time step $t$ is simulated with both the Monte-Carlo method and reachability analysis. We see that the polytope in the next time step contains all the points in the Monte-Carlo simulation, and thus the convergence of this polytope towards the equilibrium guarantees that every Monte-Carlo simulation will be stable. The main task involved in this problem is appropriately setting $A^{(t)}$ and solving for $b^{(t)}$ that can accurately predict the dynamic behavior. First, we discuss how $A^{(t)}$ can be constructed so that it maximizes the chance of converging to the equilibrium while the number of required hyperplanes stays approximately proportional to the system size.

\subsection{Template Construction}
In this section, we discuss the construction of the polytopic template $A$. The shape of this polytope is important for the convergence to equilibrium. The condition for the convergence of the polytope in a linear system is derived here using eigenvalue decomposition. Although the system is nonlinear, this construction guarantees the convergence near the equilibrium if the system is stable. The construction of the polytope is based on the eigenvalue decomposition at the equilibrium. We first linearize our system at the equilibrium, and we let the linearized system dynamic to be $\dot{x}=Jx$.
The eigenvalue decomposition in the real system representation is given by $J=Q\Lambda Q^{-1}$. We construct the polytopic template, $P\in \mathbb{R}^{m\times n}$. We denote each block of $\Lambda$ as $\Lambda_{(l)}$ so that $\Lambda=\mathrm{blkdiag}(\Lambda_{(1)},...,\Lambda_{(L)})$. For a real matrix $A$, the real system representation gives $\Lambda_{(l)}=\lambda_l$ or $\Lambda_{(l)}= \begin{bmatrix}\sigma_l & \omega_l \\ -\omega_l & \sigma_l \end{bmatrix}$, and its block diagonal can be represented as,

\begin{equation}
\Lambda
= \begin{bmatrix}
 \lambda_1 &               &               &               &             & 0 \\
                    & \ddots    &               &               &             &  \\
                    &                & \sigma_l  & \omega_l &             &  \\
                    &                & -\omega_l & \sigma_l  &             &  \\
                    &                &               &               & \ddots &  \\
                0  &                &               &               &             & \lambda_L
\end{bmatrix}.
\end{equation}

 If the associated eigenvalue is real, $\lambda_l$, then $n_l=1$ and if it is imaginary, $\sigma_l\pm j\omega_l$, the $n_l=2$.
We construct a matrix $A=\hat{A}Q^{-1}$ such that $\hat{A}=\mathrm{blkdiag}(\hat{A}_{(1)},...,\hat{A}_{(L)})$ and $\hat{A}_{(l)}\in \mathbb{R}^{m_l\times n_l}$. The construction rule is as follows:
\begin{enumerate}
\item If the eigenvalue at block $(l)$ is real, then $\hat{A}_{(l)}=\begin{bmatrix} 1 & -1 \end{bmatrix}^T$
\item If the eigenvalues at block $(l)$ are complex conjugate pair, then 
$$\hat{A}_{(l)}=
 \begin{bmatrix}
 \cos(\psi_1) & \sin(\psi_1) \\
 \vdots      & \vdots     \\
 \cos(\psi_{m_l}) & \sin(\psi_{m_l})
 \end{bmatrix} $$
where $\psi_k=\frac{2k\pi}{m_l}$ and $m_l$ is chosen to satisfy the inequality, $\tan\big(\frac{\pi}{m_l}\big) < \big|\frac{\sigma_l}{\omega_l} \big|$.
\end{enumerate}

\begin{figure}[!htbp]
	\centering
	\includegraphics[width=3.6in]{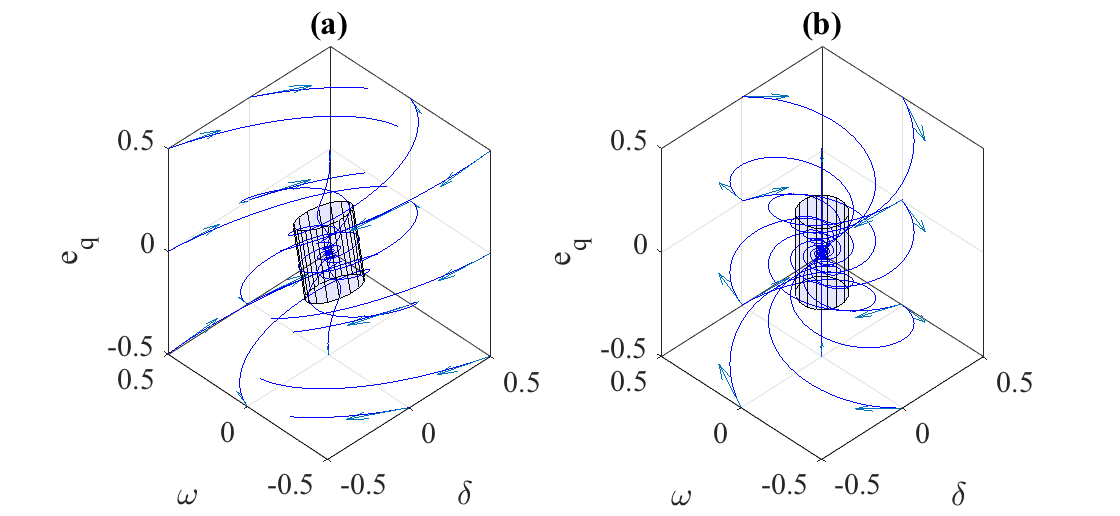}
	\caption{The polytope constructed for the 2 bus system, which will be presented in the result section. The constructed polytope as well as the original dynamics for $\dot{x}=Jx$ is shown in (a) and $\dot{\tilde{x}}=\Lambda \tilde{x}$ in (b) }
	\label{cylinder}
\end{figure}

The use of eigenvalue decomposition allows us to develop an efficient method to build template that is guaranteed to converge near the equilibrium while the number of planes is limited by $pn$ where $n$ is the system size and $p=\max_l\big[\pi/\tan^{-1}(|\sigma_l/\omega_l|)$\big]. The constructed polytope will become an intersection of cylinders and intervals along the eigenvectors of the system and this is illustrated in Figure \ref{cylinder}. Figure \ref{cylinder} (a) shows the constructed polytope in the original state space, and Figure \ref{cylinder} (b) shows the polytope in linear transformed space such that the eigenvectors aligns with the axis.

\begin{figure}[!htbp]
	\centering
	\includegraphics[width=3.6 in]{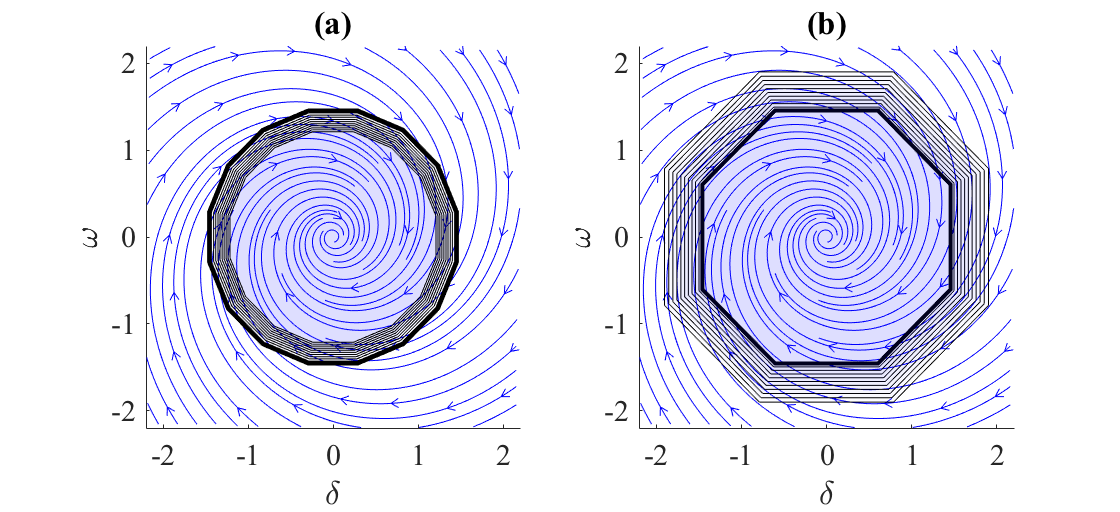}
	\caption{The result of reachability analysis when (a) the condition, $\tan\big(\frac{\pi}{m_l}\big) < \big|\frac{\sigma_l}{\omega_l} \big|$, is met and (b) the condition is not met. The polytope with the thick line is the initial polytope.}
	\label{round}
\end{figure}

The condition on the number of hyperplanes, $\tan\big(\frac{\pi}{m_l}\big) < \big|\frac{\sigma_l}{\omega_l} \big|$, along the complex eigenvalues ensures that the polytope approximates the 2-norm tightly enough so that the polytope is invariant. Figure \ref{round} shows the converging and diverging polytope depending on the number of hyperplane requirement is met or not.

\subsection{Dynamic template}
\begin{figure}[!htbp]
	\centering
	\includegraphics[width=3.2in]{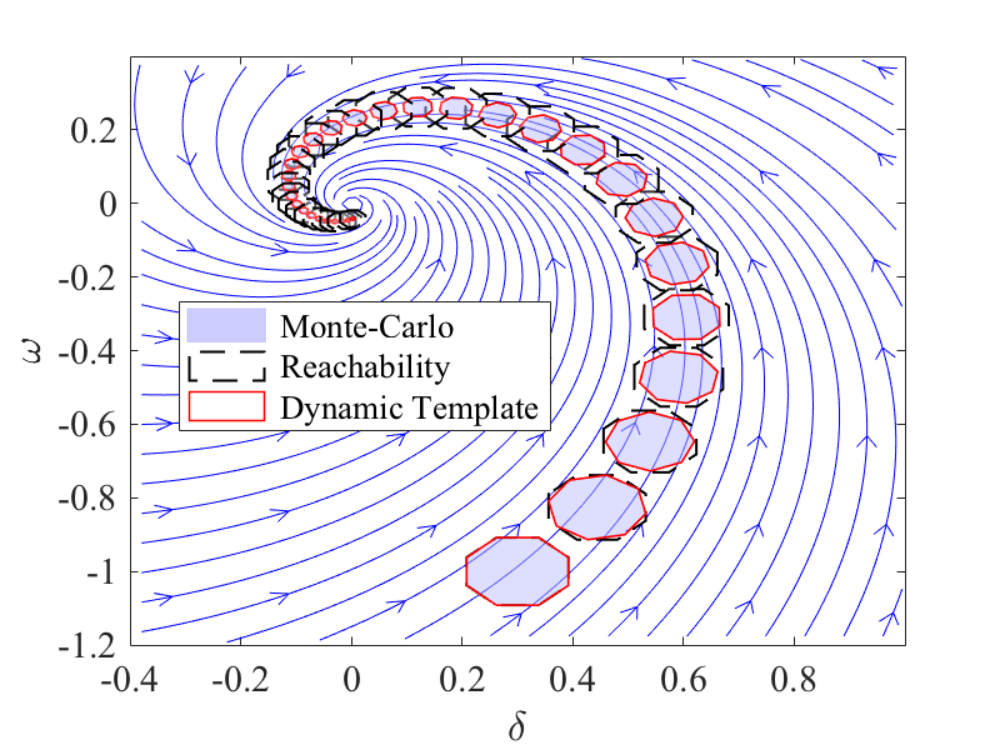}
	\caption{Reachability analysis on an under-damped linear system.}
	\label{illust_linear_full}
\end{figure}

After constructing the initial template, the template is updated at every time step to capture the change in template due to the dynamics. The dynamic template approach was introduced in \cite{sassi12}, and the update rule is given as follows:
\begin{equation}
A^{(t+1)}=A^{(t)}\cdot (I+\Delta t \cdot J(\tilde{x}^{(t)}))^{-1}
\end{equation}
where $\tilde{x}^{(t)}$ is the centroid of the polytope at time step $t$, and each row of $A$ is renormalized. Alternative to the centroid, the center point of initial polytope can be simulated and used as the point for computing the Jacobian. This update captures the linear component of the dynamics and adapt to the change in the orientation, which can significantly reduce the wrapping effect coming from the limited number of hyperplanes. In a linear system, the reachability analysis with dynamic template computes the exact states. Figure \ref{illust_linear_full} shows reachability analysis with both fixed template and dynamic template as well as the Monte-Carlo approach. The bound computed with the dynamic template is exact to the Monte-Carlo simulation as stated. The gap between the blacked dashed line and red straight line in Figure \ref{illust_linear_full} is the wrapping effect caused by enforcing the fixed polytope. This effect is resolved with the dynamic template.

While the dynamic template produce tighter bound than fixed template, it may not be guaranteed to converge to the equilibrium after going through the nonlinearities. To address this concern, we concatenate the original constructed template and its duplicate, and update the duplicate template with the dynamic template. In the next section, we formulate an optimization problem to solve the exact polytope at every time step by solving $b^{(t)}$.

\subsection{Problem Formulation}
After building the fixed polytopic template $A$ and applying the time stepping technique to the dynamics, the convergence of the polytope is tracked by computing the bound, $b$. We consider the following optimization problem in order to bound the reachable half space and apply it to every hyperplanes in the polytope:

\begin{equation}
	\begin{aligned}
		\rho_i^{(t+1)} = \underset{x,y}{\text{max}} \hskip 1em & A_i\cdot x^{(t+1)} \\
		\text{subject to} \hskip 1em & x^{(t+1)}=x^{(t)} +\Delta t f(x^{(t)},y^{(t)}) \\
		&g(x^{(t)},y^{(t)}) =0 \\
		& Ax^{(t)} \leq b^{(t)}.
	\end{aligned}
	\label{opt_original}
\end{equation}
This formulation computes the tightest polytope that contain all possible reachable space subject to the system dynamics and current state bounds. This is a non-convex problem and obtaining the exact solution is intractable for a large scale system. We assume that the nonlinearity is contained within the differential variables, and thus the function $f(x,y)$ and $g(x,y)$ can be written in the following form:

\begin{equation}
	\begin{aligned}
		f(x^{(t)},y^{(t)})&=\tilde{f}(x^{(t)},y^{(t)})+f_xx^{(t)}+f_yy^{(t)} \\
		g(x^{(t)},y^{(t)})&=\tilde{g}(x^{(t)})+g_xx^{(t)}+g_yy^{(t)}.
	\end{aligned}
\end{equation}
where $\tilde{f}(x^{(t)},y^{(t)})$ and $\tilde{g}(x^{(t)})$ contains the nonlinear terms in differential and algebraic equations respectively. We introduce new variables for the nonlinear terms

\begin{equation}
	\begin{aligned}
		u&=\tilde{f}(x^{(t)},y^{(t)}) \\
		v&=\tilde{g}(x^{(t)})
	\end{aligned}
\end{equation}
Then the original optimization problem \ref{opt_original} becomes

\begin{equation}
	\begin{aligned}
		\underset{x,y,u,v}{\text{maximize}} \hskip 1em & A_i\cdot x^{(t+1)} \\
		\text{subject to} \hskip 1em & x^{(t+1)}=x^{(t)} +\Delta t (u+f_xx^{(t)}+f_yy^{(t)}) \\
		& v+g_xx^{(t)}+g_yy^{(t)} =0 \\
		& Ax^{(t)} \leq b^{(t)} \\
		& u=\tilde{f}(x^{(t)},y^{(t)}) \\
		& v=\tilde{g}(x^{(t)})
	\end{aligned}
	\label{opt_original2}
\end{equation}

In this formulation, the nonlinear terms are confined in the variables $u$ and $v$ where the relaxation can be applied to solve the optimization problem efficiently. 

\section{LP Relaxation}
To take advantage of scalability and reliability of linear programming, the nonlinear constraints in the formulation needs to be replaced with linear constraints. These linear constraints are used as the outer-approximation of the nonlinear functions. In this section, we present affine envelopes for bilinear and sinusoidal functions, which can efficiently bound the non-linearities.

\subsection{McCormick Envelopes}
McCormick envelop is an outer-approximation for a bilinear function, $u=xy$, given its bounds, $\underbar{$x$}\leq x\leq\overline{x}$ and $\underbar{$y$}\leq y\leq\overline{y}$ \cite{mccormick76}. The following linear constraints replace the nonlinear function:

\begin{equation}
	\begin{aligned}
		u & \geq \underbar{$x$}y+x\underbar{$y$}-\underbar{$x$}\underbar{$y$} \\
		u & \geq \overline{x}y+x\overline{y}-\overline{x}\overline{y} \\
		u & \leq \overline{x}y+x\underbar{$y$}-\overline{x}\underbar{$y$} \\
		u & \leq x\overline{y}+\underbar{$x$}y-\underbar{$x$}\overline{y}.
	\end{aligned}
	\label{eqn_McCormick}
\end{equation}

\subsection{Sinusoidal Envelopes}
Similar to McCormick envelopes, the outer-approximation is developed in this paper for sinusoidal functions to replace the nonlinear function by linear inequalities. Given any sinusoidal function $u=f(\delta)=\sin(\delta+\phi)$ with any phase shift, $\phi$, we define the slope of the chord between points $a$ and $b$ as $m_{a,b}=\big(\frac{f(a)-f(b)}{a-b}\big)$. Given $\underline{\delta} \leq \delta \leq \overline{\delta}$ and $\overline{\delta}-\underline{\delta} \leq \frac{\pi}{2}$, the linear envelope includes $|u| \leq 1$ and add additional inequalities based on the following two cases.

\subsubsection{Convex/Concave Region}

If $f'(\overline{\delta}) \geq m_{\overline{\delta},\underline{\delta}} $ and $f'(\underline{\delta}) \leq m_{\overline{\delta},\underline{\delta}}$, the function is convex in the given region. In this case, the outer-approximation is built as follows:

\begin{equation}
	\begin{aligned}
		u & \geq f'(\overline{\delta})(\delta-\overline{\delta})+f(\overline{\delta})  \\
		u & \geq f'(\underline{\delta})(\delta-\underbar{$\delta$})+f(\underline{\delta}) \\
		u & \leq m_{\overline{\delta},\underline{\delta}} (\delta-\overline{\delta})+f(\overline{\delta}) \\
		u & \geq m_{\overline{\delta},\underline{\delta}} (\delta-f'^{-1}(m_{\overline{\delta},\underline{\delta}}))+f(f'^{-1}(m_{\overline{\delta}.\underline{\delta}}))
	\end{aligned}
	\label{eqn_sinenv_convex}
\end{equation}
The first inequality is from the definition of convex function and the second inequality is from the mean value theorem. The third and fourth inequalities are first order condition applied at the boundary points. If $f'(\overline{\delta}) \leq m_{\overline{\delta},\underline{\delta}} $ and $f'(\underline{\delta}) \geq m_{\overline{\delta},\underline{\delta}}$, then the function is concave, and the signs of inequalities in Equation \ref{eqn_sinenv_convex} flip to the other side.

\subsubsection{Monotonic Region}
The function is monotonically increasing if $f'(\overline{\delta}) \leq m_{\overline{\delta},\underline{\delta}} $ or $f'(\underline{\delta}) \leq m_{\overline{\delta},\underline{\delta}}$, and the outer-approximation of this region can be written as
\begin{equation}
	\begin{aligned}
		u & \geq f'(\overline{\delta})(\delta-\overline{\delta})+f(\overline{\delta})  \\
		u & \leq f'(\underline{\delta})(\delta-\underbar{$\delta$})+f(\underline{\delta})  \\
		u & \leq m_{\overline{\delta},\overline{\eta}}(\delta-\overline{\delta})+f(\overline{\delta}) \\
		u & \geq m_{\underline{\eta},\underline{\delta}}(\delta-\underline{\delta})+f(\underline{\delta}).  \\
	\end{aligned}
	\label{eqn_sinenv_monotone}
\end{equation}
Similar to convex/concave relationship, the function is monotonically decreasing if $f'(\overline{\delta}) \geq m_{\overline{\delta},\underline{\delta}} $ or $f'(\underline{\delta}) \geq m_{\overline{\delta},\underline{\delta}}$, and the signs of inequalities in Equation \ref{eqn_sinenv_monotone} flip to the other side.

\begin{figure}[!htbp]
	\centering
	\includegraphics[width=3.6in]{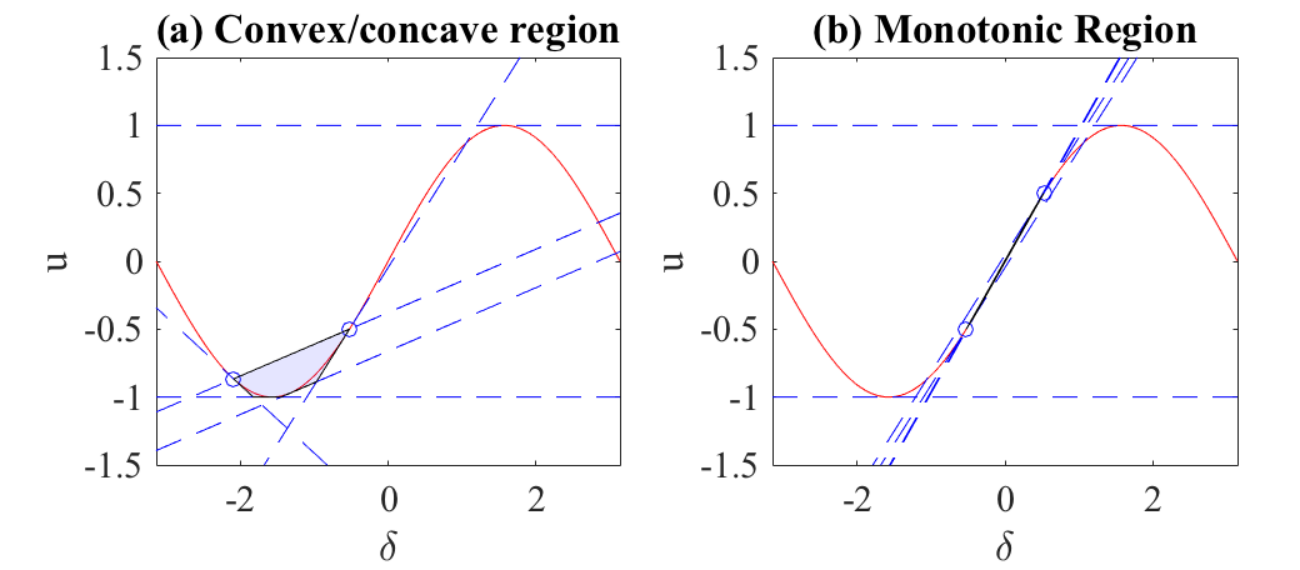}
	\caption{outer-approximation for sinusoidal function. The given bound is marked with blue circle, and the polytope approximation is colored with light blue.}
	\label{fig_sin_env}
\end{figure}

Figure \ref{fig_sin_env} shows examples of the two cases. We only consider tight input bound, $\overline{\delta}-\underline{\delta} \leq \frac{\pi}{2}$, since large bound causes high gap between the original and relaxed problems. The tightness of the input bound is set as the exit criteria to conclude instability or inability to conclude stability of the system. Now, the set of linear inequalities form a closed set that contains the nonlinear function, and we denote this set by $\Gamma$ where

\begin{equation}
	\Gamma(\underline{x},\overline{x},\underline{y},\overline{y},x,y)=\{u| \text{ Equation (\ref{eqn_McCormick}, \ref{eqn_sinenv_convex}, \ref{eqn_sinenv_monotone})}\}.
	\label{eqn_Gamma}
\end{equation}
This outer-approximation will take the bounds of the input and form a linear inequalities as a function of the inputs and nonlinear outputs.

\subsection{Bounds on nonlinear functions}
The outer-approximation requires the bounds on both differential and algebraic variables, $x^{(t)}$ and $y^{(t)}$. Since the polytope $\Omega^{(t)}$ is closed and the algebraic equations is linear with respect to the algebraic variables, the upper bound and under bound can be computed prior with the following optimization problems. The upper-bound on the differential variables can be simply computed by

\begin{equation}
	\begin{aligned}
		\overline{x}_i=\underset{x}{\text{max}} \hskip 1em & e_i\cdot x^{(t)} \\
		\text{subject to} \hskip 1em & Ax^{(t)} \leq b^{(t)}.
	\end{aligned}
	\label{opt_x_max}
\end{equation}

The under-bound of differential variable, $\underbar{$x$}$, can be computed by solving minimization problem of Equation \ref{opt_x_max}. The upper-bound on the algebraic variables can be computed by

\begin{equation}
	\begin{aligned}
		\overline{y}_i=\underset{x,y,v}{\text{max}} \hskip 1em & e_i\cdot y^{(t)} \\
		\text{subject to} \hskip 1em & v+g_xx^{(t)}+g_yy^{(t)} =0 \\
		& Ax^{(t)} \leq b^{(t)} \\
		& v\in\Gamma_g(\underbar{$x$}^{(t)},\overline{x}^{(t)},x^{(t)}).
	\end{aligned}
	\label{opt_y_max}
\end{equation}
Similarly, the under bound of differential variable, $\underbar{$y$}$, can be computed by solving minimization problem of Equation \ref{opt_y_max}. These problems should be solved in the beginning of each time step, and the solution can be reused within the same time step.

\subsection{LP Relaxation}
Based on the outer-approximation proposed in the previous section, we replace the nonlinear terms in optimization problem in Equation \ref{opt_original2} with linear inequality constraints,

\begin{equation}
	\begin{aligned}
		\underset{x,y,u,v}{\text{maximize}} \hskip 1em & A_i\cdot x^{(t+1)} \\
		\text{subject to} \hskip 1em & x^{(t+1)}=x^{(t)} +\Delta t (u+f_xx^{(t)}+f_yy^{(t)}) \\
		& v+g_xx^{(t)}+g_yy^{(t)} =0 \\
		& Ax^{(t)} \leq b^{(t)} \\
		& u\in\Gamma_f(\underbar{$x$}^{(t)},\overline{x}^{(t)},\underbar{$y$}^{(t)},\overline{y}^{(t)},x^{(t)},y^{(t)}) \\ 
		& v\in\Gamma_g(\underbar{$x$}^{(t)},\overline{x}^{(t)},x^{(t)})
	\end{aligned}
	\label{opt_relaxed}
\end{equation}
which can be solved using linear programming.
Suppose the solution for the optimization problem in Equation \ref{opt_relaxed} is $b_i^{(t+1)}$. Since the linear envelop is a relaxation of the original problem,
\begin{equation}
	Ax^{(t+1)}\leq \rho^{(t+1)} \leq b^{(t+1)}
\end{equation}
Thus,
\begin{equation}
	\Omega_{(A,\rho)}^{(t+1)} \subseteq \Omega_{(A,b)}^{(t+1)}
\end{equation}

Therefore, the computed polytope with relaxed problem is guaranteed to contain the polytope from solving the original problem in Equation \ref{opt_original}.

\subsection{Criteria for stability}
The exit condition is defined in this section to determine whether the states converged to the equilibrium or not. The simulation could be terminated if all the states are close to the equilibrium by an arbitrary value $\epsilon$, and this condition for the convergence is $Ax-b^{(t)}<\epsilon$. An alternative approach would be computing $\epsilon$ that will certify the convergence in all states if they satisfy $b^{(t)}<\epsilon b^{(0)}$. This convergence can be certified by forming computing the largest invariant set. The invariance of the polytope can be determined by solving the following optimization problem,

\begin{equation}
	\begin{aligned}
		\mu_i^\epsilon=\underset{x,y,u,v}{\text{min}} \hskip 1em & A_i\cdot f(x,y) \\
		\text{subject to} \hskip 1em & f(x,y)=u+f_xx+f_yy \\
		& v+g_xx+g_y=0 \\
		& Ax \leq \epsilon b^{(0)} \\
		& A_ix = \epsilon b_i^{(0)} \\
		& u\in\Gamma_f(\underbar{$x$},\overline{x},\underbar{$y$},\overline{y},x,y) \\ 
		& v\in\Gamma_g(\underbar{$x$},\overline{x},x).
	\end{aligned}
	\label{opt_exit_criteria}
\end{equation}

If $\mu^\epsilon=\max_i\mu_i^\epsilon$ is negative, then it certifies that $A(x-x_{eq})\leq \epsilon b^{(0)}$ is invariant. If we define $V=max_i A_i(x-x_{eq})$ as a level set, then $\mu_\epsilon=\max_{x\in \partial\mathcal{P}^\epsilon}\frac{\partial V}{\partial x}f(x)$ where $\partial\mathcal{P}^\epsilon$ is the boundary of $A\leq \epsilon b$. The condition $\mu^\epsilon<0$ indicate that all the dynamic has direction inward to the polytope, and the invariance of the set, $A\leq \epsilon b$, can be concluded. A bisection algorithm can be applied to find the maximum $\epsilon$ where $\mu_i^\epsilon$ can be used as the oracle to determine the invariance. To set the exit condition for case where the polytope does not converge, we use the assumption in our outer-approximation, which was $\delta_{max}-\delta_{min}\leq\frac{\pi}{2}$. Once this condition is violated, we exit the algorithm and declare it was unable to certify stability. Once the stability criterias are established, our final algorithm is presented below:

\begin{figure}[!htbp]
	\centering
	\includegraphics[width=2.6in]{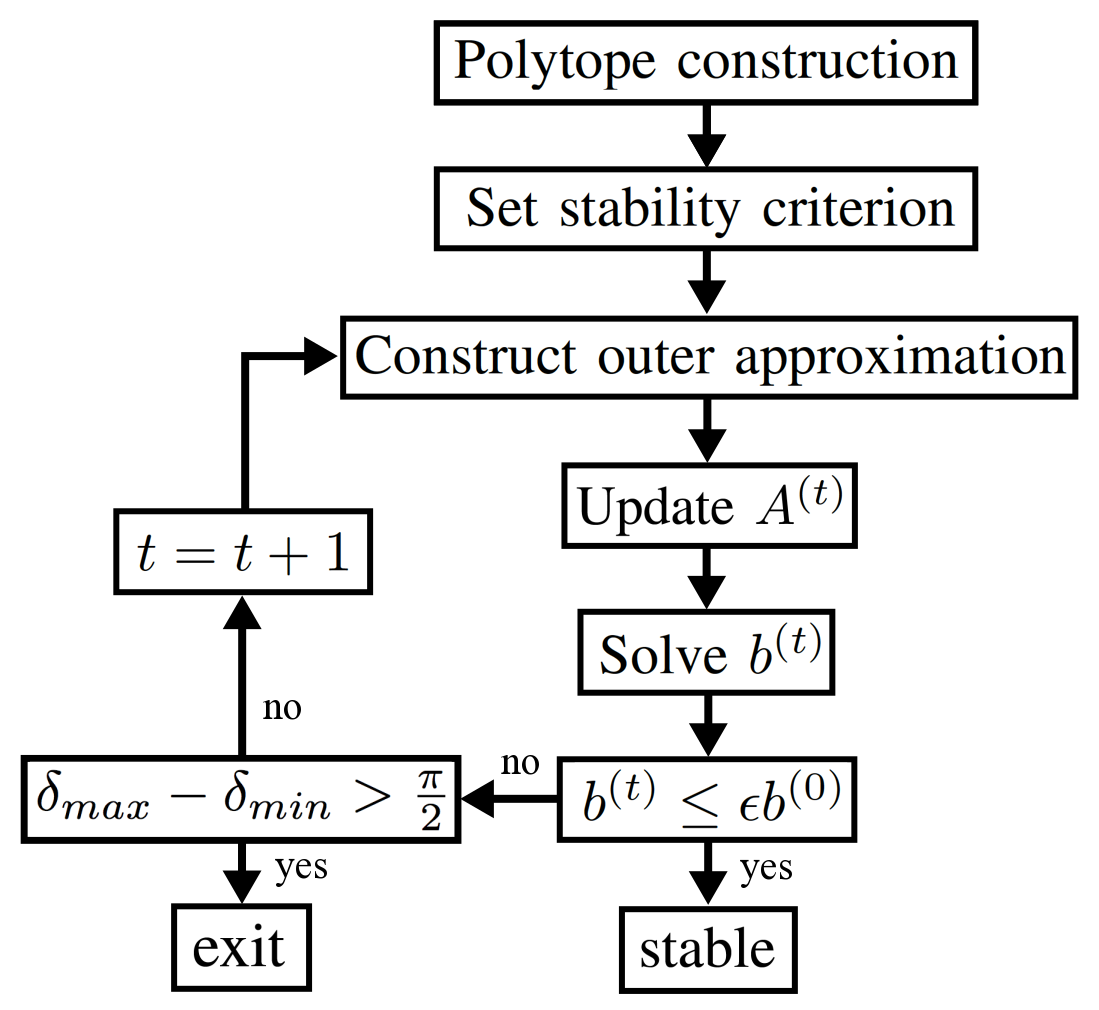}
	\caption{Flow chart for the proposed algorithm.}
	\label{overall_alg}
\end{figure}

\section{Results}
\subsection{System Model}
In this paper, we consider the third order synchronous dynamics model given in the following equation \cite{machowski11}:

\begin{equation}
	\begin{aligned}
		\dot{\delta}_i&=\omega_i \\
		\dot{\omega}_i&=\frac{1}{M_i}\big[-D_i\omega_i-P_{e,i}+P_{m,i}\big] \\
		\dot{e}_{q,i}'&=\frac{1}{T_{d,i}'}\big[E_{f,i}-e_{q,i}'-(x_{d,i}-x_{d,i}')\tilde{i}_{d,i} \big] \\
		P_{e,i}&=\tilde{v}_{d,i}\tilde{i}_{d,i}+\tilde{v}_{q,i}\tilde{i}_{q,i} \\
		\begin{bmatrix} \tilde{v}_{d,i} \\ \tilde{v}_{q,i} \end{bmatrix}&=\begin{bmatrix} 0 \\ e_{q,i}' \end{bmatrix}-\begin{bmatrix} 0 & x_{d,i}' \\ -x_{d,i}' & 0 \end{bmatrix} \begin{bmatrix} \tilde{i}_{d,i} \\ \tilde{i}_{q,i} \end{bmatrix}
	\end{aligned}
\end{equation}
where $M_i$ and $D_i$ are inertia constant, damping constant of the generator $i$, and $\tilde{v}$ and $\tilde{i}$ are the voltage and current in the reference frame with respect to the generator bus. The network is modeled in \textit{dq}-frame with a global reference:
\begin{equation}
	\begin{bmatrix} i_{d} \\ i_{q} \end{bmatrix}=\begin{bmatrix} G & -B \\ B & G \end{bmatrix} \begin{bmatrix} v_{d} \\ v_{q} \end{bmatrix}.
\end{equation}

We model the load here as a fixed current sink:
\begin{equation}
	\begin{bmatrix} i_{d,k} \\ i_{q,k} \end{bmatrix}=\begin{bmatrix} i_{d,k}^{load} \\ i_{q,k}^{load} \end{bmatrix}.
\end{equation}

The generator voltage and current in its local angle reference can be transformed into a global reference with the following relationship:
\begin{equation}
	\begin{bmatrix} \tilde{v}_{d,i} \\ \tilde{v}_{q,i} \end{bmatrix}=\begin{bmatrix} v_{i}\sin(\delta_i-\delta_i) \\ v_{i}\cos(\delta_i-\delta_i) \end{bmatrix}=\begin{bmatrix} \sin\delta_i & \cos\delta_i \\ -\cos\delta_i & \sin\delta_i \end{bmatrix} \begin{bmatrix} v_{d,i} \\ v_{q,i} \end{bmatrix}.
\end{equation}

We can replace the voltage and current in local generator angle reference by the global reference as follows:
\begin{equation}
	\begin{aligned}
		\dot{e}_{q,i}'&=\frac{1}{T_{d,i}'}\big[E_{f,i}-e_{q,i}' \\
		& \hskip 2em -(x_{d,i}-x_{d,i}')(i_{d,i}\cos\delta_i+i_{q,i}\sin\delta_i) \big] \\
		P_{e,i}&=v_{d,i}i_{d,i}+v_{q,i}i_{q,i}  \\
		\begin{bmatrix} v_{d,i} \\ v_{q,i} \end{bmatrix}&=\begin{bmatrix} \sin\delta_i & -\cos\delta_i \\ \cos\delta_i & \sin\delta_i \end{bmatrix} \begin{bmatrix} 0 \\ e_{q,i}' \end{bmatrix} \\
		& \hskip 8em -\begin{bmatrix} 0 & x_{d,i}' \\ -x_{d,i}' & 0 \end{bmatrix} \begin{bmatrix} i_{d,i} \\ i_{q,i} \end{bmatrix}.
	\end{aligned}
\end{equation}

This transformation enables writing the algebraic equation to be linear with respect to the algebraic variable.

\subsection{2 bus system}
\begin{figure}[!htbp]
	\centering
	\includegraphics[width=3.2in]{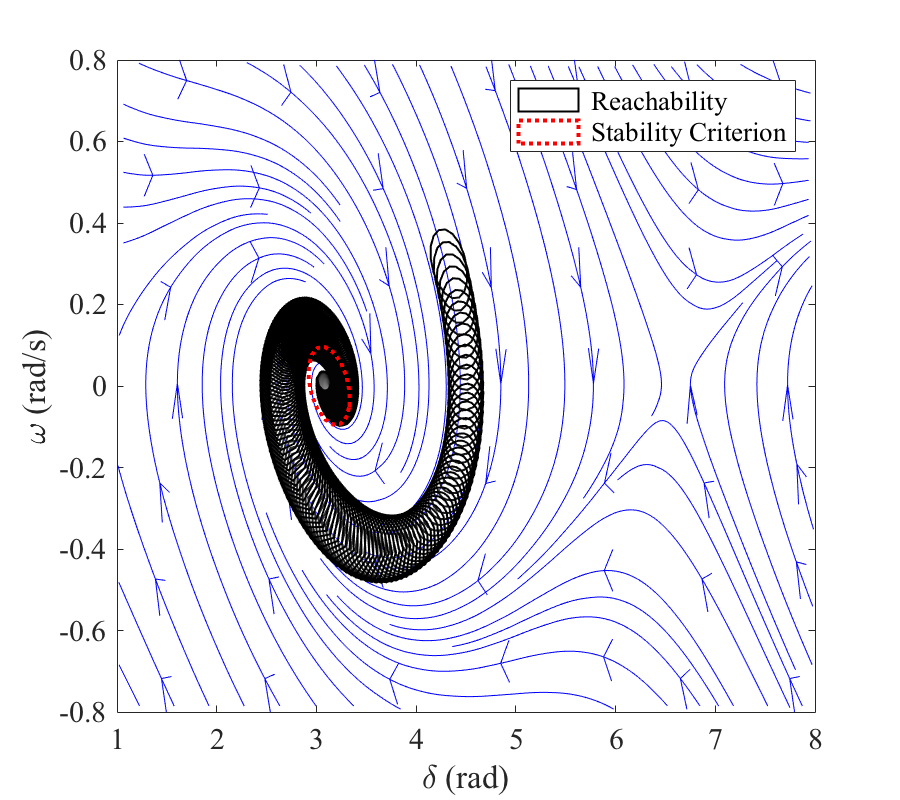}
	\caption{The reachability analysis on 2 bus system and the phase portrait of the system is shown.}
	\label{phase_portrait_stable}
\end{figure}

The proposed simulation was performed on a 2 bus system with a single generator to illustrate and visualize our approach. The linearization gave the eigenvalues of the system at the equilibrium to be $\lambda=\{-0.153\pm j0.433,-0.128\}$, and the constructed initial polytope was illustrated in Figure \ref{cylinder} (a). Figure \ref{phase_portrait_stable} shows the result of our simulation on top of the phase portrait of the system. We determine the stability criteria coefficient and marked with a red dotted line where we terminate the analysis once the polytope gets inside the stopping criteria.

\begin{figure}[!htbp]
	\centering
	\includegraphics[width=3.2in]{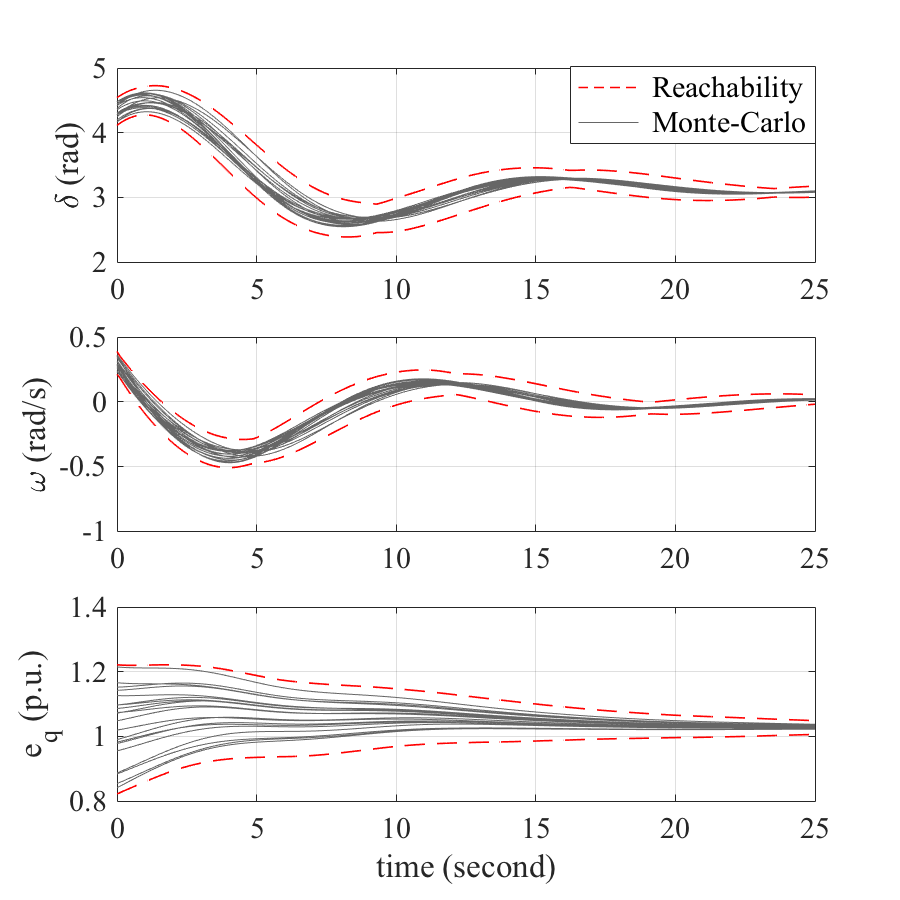}
	\caption{The Monte-Carlo simulation of random initial condition within the initial polytope is presented as well as the bound given by the reachability analysis.}
	\label{time_domain_stable}
\end{figure}

In Figure \ref{time_domain_stable}, the Monte-Carlo simulation as well as the bound from the reachability analysis is shown. All of the simulations stays within this bound as long as the initial value was within the initial polytope of the reachability analysis.

\begin{figure}[!htbp]
	\centering
	\includegraphics[width=3.2in]{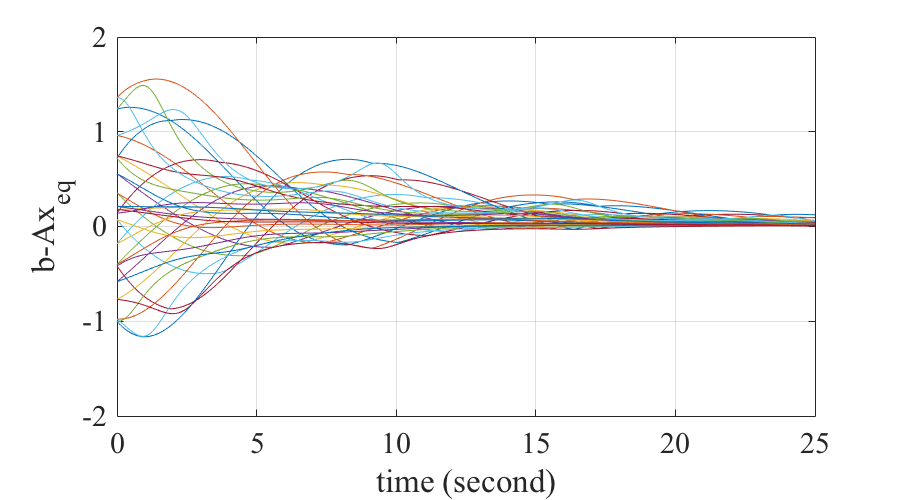}
	\caption{The evolution of the polytope is illustrated with $b-Ax_{eq}$. If every simulation points converge to the equilibrium, these values should converge to zero as well.}
	\label{b_stable}
\end{figure}

Figure \ref{b_stable} shows the distance of hyperplanes in the polytope from the stable equilibrium. We see that all the hyperplanes eventually converge to the equilibrium

\subsection{14 bus system}
In this section, we consider IEEE 14 bus and 39 bus systems to demonstrate our approach in larger systems. We relocate the current point to be away from the equilibrium and construct a polytope around the initial state. This approach generalizes the simulation of a line or generator contingency where the stable equilibrium changes to a different point. The optimization problem was solved with YALMIP in MATLAB \cite{lofberg04}. For 14 bus system, the reachability analysis was performed with the time step of 100 $ms$ for the duration of 25 seconds. The total solver time in YALMIP was 177.8 seconds where 137.3 seconds were used to solve the optimization problem in Equation \ref{opt_relaxed}, and 40.5 seconds were used to construct the outer-approximation. Computation of the stopping criteria took 29.64 seconds in YALMIP solver time.

\begin{figure}[!htbp]
	\centering
	\includegraphics[width=3.2in]{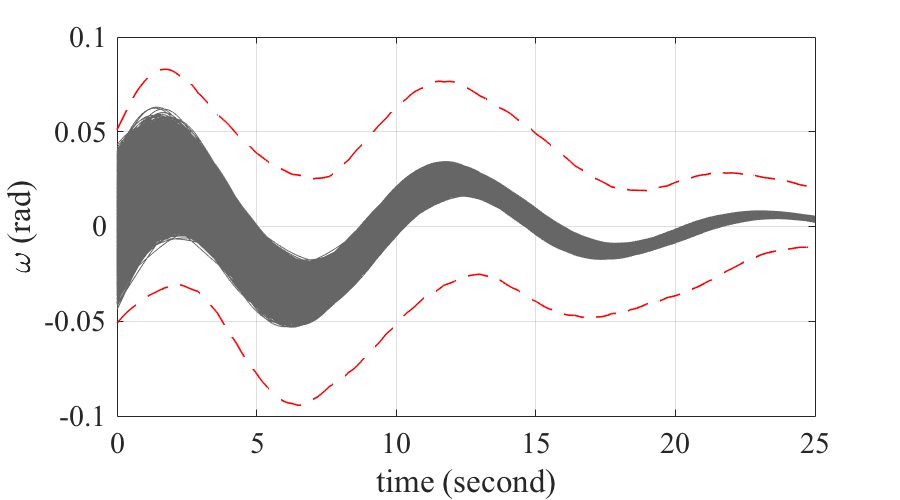}
	\caption{The rotor angle frequency at bus 3 where the grey lines are samples from the Monte-Carlo simulation and red is the bound computed by the reachability analysis for 14 bus system.}
	\label{state_14}
\end{figure}

\begin{figure}[!htbp]
	\centering
	\includegraphics[width=3.2in]{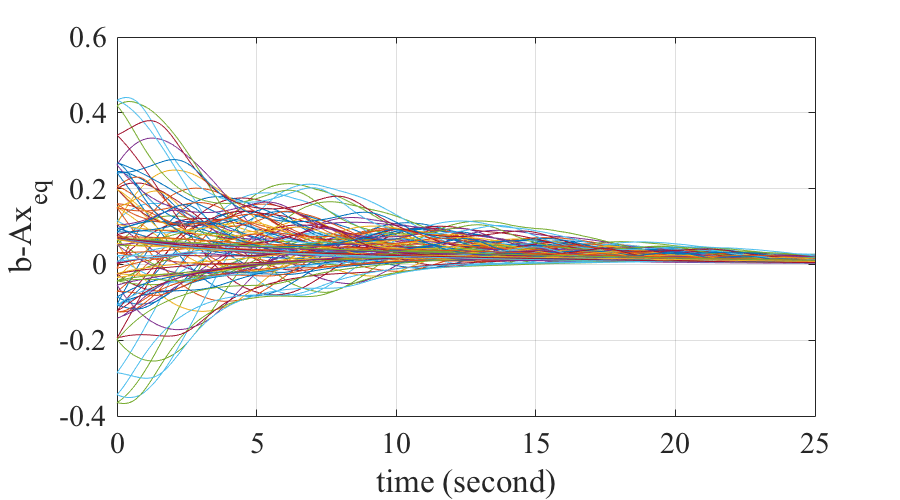}
	\caption{Distance of hyperplanes from the stable equilibrium in 14 bus system.}
	\label{b_14}
\end{figure}

\begin{figure}[!htbp]
	\centering
	\includegraphics[width=3.2in]{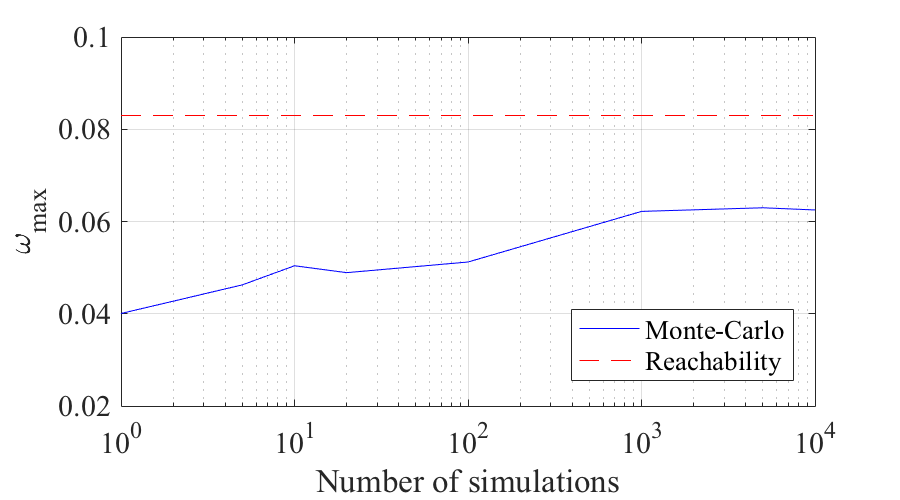}
	\caption{Comparison between maximum rotor frequency estimation with the Monte-Carlo simulation and reachability analysis for 14 bus system.}
	\label{result_comparison_14}
\end{figure}

While the convergence of reachability analysis guarantees stability against uncertain initial conditions, the analysis provides additional information about the trajectories. While sampling-based approaches give estimation of the maximum or minimum of a state along the trajectory, every solution is under-estimated. The naive Monte-Carlo simulation was implemented where the trajectory was computed with random initial condition drawn from the initial polytope with uniform distribution. The Monte-Carlo simulation takes about 0.04 seconds per simulation, and the result is compared in Figure \ref{result_comparison_14}. The advantage of reachability analysis is that it can bound the maximum and minimum while the Monte-Carlo approach requires a large number of simulations to find an accurate solution.

\subsection{39 bus system}
For 39 bus system with 15 seconds and the time step of 100 $ms$, the total solver time was 423.5 seconds where 312.5 seconds were used to solve the optimization problem in Equation \ref{opt_relaxed} and 111 seconds were used to construct the outer-approximation. 
The stopping criteria took 230.7 seconds where the epsilon was computed to be 0.1313.

\begin{figure}[!htbp]
	\centering
	\includegraphics[width=3.2in]{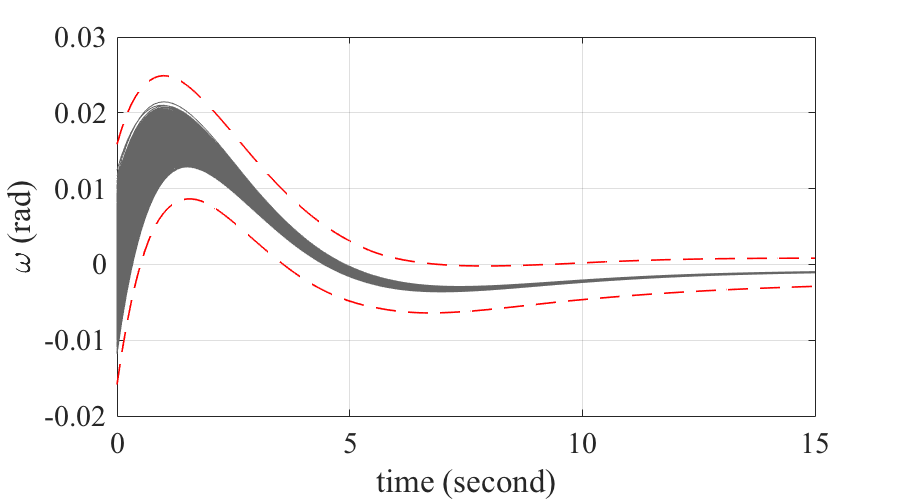}
	\caption{The rotor angle frequency at bus 39 with the Monte-Carlo simulation (grey) and reachability analysis (red) for 39 bus system.}
	\label{state_39}
\end{figure}

Similar to 14 bus system case study, the estimation of maximum bound on the rotor angle was compared between the Monte-Carlo simulation and reachability analysis. The Monte-Carlo simulation takes about 0.09 seconds per simulation in 39 bus system.

\begin{figure}[!htbp]
	\centering
	\includegraphics[width=3.2in]{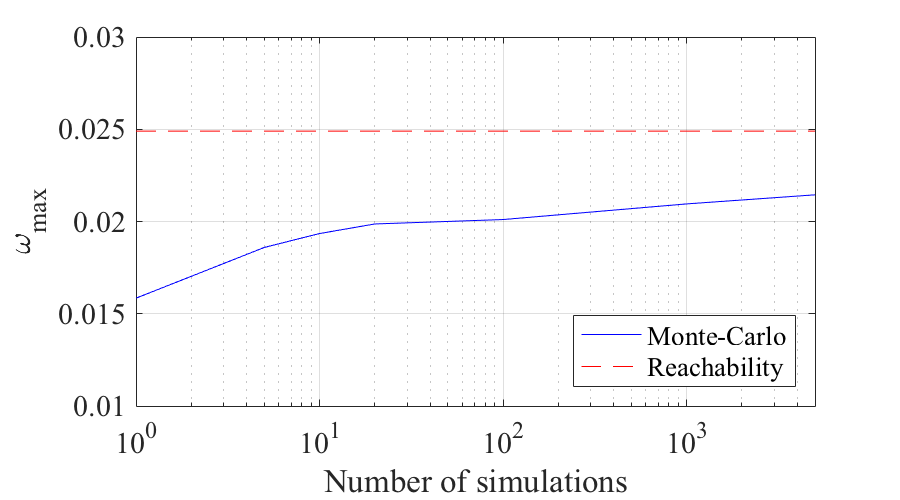}
	\caption{Comparison between maximum rotor frequency estimation with Monte-Carlo simulation and reachability analysis for 39 bus system.}
	\label{result_comparison_39}
\end{figure}

\section{Discussion}
Stability assessment based on reachability analysis is enabled from the development of linear programming as a reliable and scalable tool. This technique allows certification of the convergence of trajectories from a neighborhood of initial states. A side-product of the convergence is the retrieval of the maximum and minimum of the states along the trajectories, which can identify violations of system constraints during transient events.

To take advantage of linear programming, the polytope was chosen to be a template that describes the cloud of states at a given time. Based on the eigenvalue decomposition, the number of hyperplanes in the polytope was approximately proportional to the size of the system while guaranteeing convergence near the equilibrium.

While our approach allows incorporation of complex generator models, one of the limitations comes from the accumulation of errors throughout the simulation. There is an inherent gap between the actual cloud of states and the computed polytope states due to solving the relaxed problem and enforcing the polytopic template. This gap limits our approach from the convergence in the case where the cloud gets too large due to trajectories diverging significantly or having large initial uncertainties.

We note that the computation of $b$ in the algorithm flow chart in Figure \ref{overall_alg} does not have to be computed sequentially. The parallel computing implementation with greater computing capabilities can make this step very efficient.

It is also possible to obtain estimates of minimum and maximum bus voltages with increased computational costs. While it is easy to infer the maximum and minimum from differential variables since the states are directly bounded by the polytope, the algebraic variables require additional computations to estimate the bounds through the algebraic equations.

\section{Conclusions}
In this paper, we presented the reachability analysis approach for transient stability assessment in power systems. The linear programming relaxation based on the outer-approximation allows this method to be more scalable and tractable than existing approaches in reachability analysis. Through this technique, we are able to certify the convergence of the system under uncertain initial conditions and obtain the bounds on the trajectories. 

An interesting future direction is the consideration of the uncertain parameters in the system. The current formulation allows incorporation of parameter uncertainty quite easily, and there are few barriers to incorporating uncertainties as random variables. In addition, the efficiency of the technique could be improved by considering Euler's implicit method instead of Euler's explicit method. This modification could allow increased step size with potentially better convergence properties.





\ifCLASSOPTIONcaptionsoff
\newpage
\fi



%

\bibliographystyle{IEEEtran}
\bibliography{references}
\nocite{*} 

%

\begin{IEEEbiography}[{\includegraphics[width=1in,height=1.25in,clip,keepaspectratio]{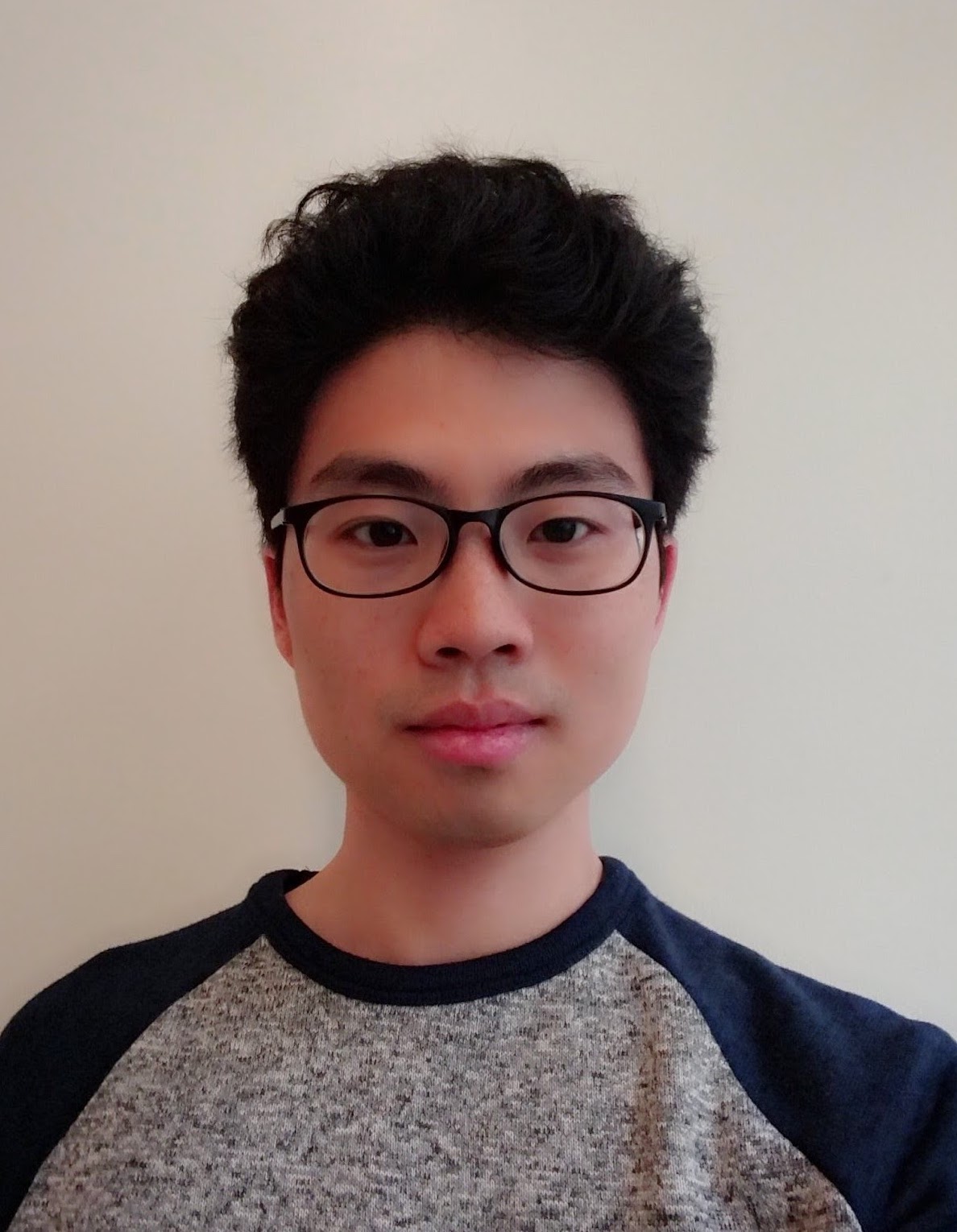}}]{Dongchan Lee}
received the M.A.Sc. degree in Electrical and Computer Engineering and B.A.Sc in Engineering Science from the University of Toronto, Canada, in 2016 and 2014, respectively. He is currently pursuing his Ph.D. degree in the Department of Mechanical Engineering at Massachusetts Institute of Technology (MIT), Cambridge. He spent a year as a power systems analyst intern at the Independent Electricity System Operator (IESO), Ontario, in 2013. His research has received the best paper award at the 2015 IEEE Electrical Power and Energy Conference. His current research interests include power systems operation, control and stability assessment. 
\end{IEEEbiography}

\begin{IEEEbiography}[{\includegraphics[width=1in,height=1.25in,clip,keepaspectratio]{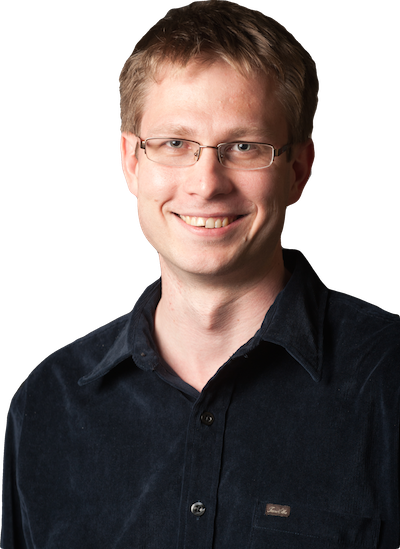}}]{Konstatin Turitsyn}
received the M.Sc. degree in physics from Moscow Institute of Physics and Technology and the Ph.D. degree in physics from Landau Institute for Theoretical Physics, Moscow, in 2007. Currently, he is an Associate Professor at the Mechanical Engineering Department of Massachusetts Institute of Technology (MIT), Cambridge. Before joining MIT, he held the position of Oppenheimer fellow at Los Alamos National Laboratory, and Kadanoff-Rice Postdoctoral Scholar at University of Chicago. His research interests encompass a broad range of problems involving nonlinear and stochastic dynamics of complex systems. Specific interests in energy related fields include stability and security assessment, as well as integration of distributed and renewable generation.
\end{IEEEbiography}






\end{document}